
\magnification1200
\input amstex.tex
\documentstyle{amsppt}

\hsize=12.5cm
\vsize=18cm
\hoffset=1cm
\voffset=2cm

\footline={\hss{\vbox to 2cm{\vfil\hbox{\rm\folio}}}\hss}
\nopagenumbers
\def\DJ{\leavevmode\setbox0=\hbox{D}\kern0pt\rlap
{\kern.04em\raise.188\ht0\hbox{-}}D}

\def\txt#1{{\textstyle{#1}}}
\baselineskip=13pt
\def\hf{{\textstyle{1\over2}}}
\def\a{\alpha}
\def\d{{\,\roman d}}
\def\e{\varepsilon}\def\E{{\roman e}}

\def\G{\Gamma}

\def\={\;=\;}

\def\zt{\zeta(\hf+it)}

\def\D{\Delta}
\def\E{{\roman e}}

\def\z{\zeta}

\def\hf{{\textstyle{1\over2}}}
\def\txt#1{{\textstyle{#1}}}

\def\le{\leqslant} \def\ge{\geqslant}
\font\tenmsb=msbm10
\font\sevenmsb=msbm7
\font\fivemsb=msbm5
\newfam\msbfam
\textfont\msbfam=\tenmsb
\scriptfont\msbfam=\sevenmsb
\scriptscriptfont\msbfam=\fivemsb

\font\ff=cmr8
\def\txt#1{{\textstyle{#1}}}
\baselineskip=13pt

\font\teneufm=eufm10
\font\seveneufm=eufm7
\font\fiveeufm=eufm5
\newfam\eufmfam
\textfont\eufmfam=\teneufm
\scriptfont\eufmfam=\seveneufm
\scriptscriptfont\eufmfam=\fiveeufm
\def\mathfrak#1{{\fam\eufmfam\relax#1}}

\font\tenmsb=msbm10
\font\sevenmsb=msbm7
\font\fivemsb=msbm5
\newfam\msbfam
     \textfont\msbfam=\tenmsb
      \scriptfont\msbfam=\sevenmsb
      \scriptscriptfont\msbfam=\fivemsb

  \def\rightheadline{{\hfil{\ff
  On some results for $|\zt|$ and a divisor problem}\hfil\tenrm\folio}}

  \def\leftheadline{{\tenrm\folio\hfil{\ff
   Aleksandar Ivi\'c }\hfil}}
  \def\emptyheadline{\hfil}
  \headline{\ifnum\pageno=1 \emptyheadline\else
  \ifodd\pageno \rightheadline \else \leftheadline\fi\fi}

\topmatter
\title
ON SOME MEAN VALUE RESULTS FOR THE ZETA-FUNCTION AND A DIVISOR PROBLEM
\endtitle
\author   Aleksandar Ivi\'c  \endauthor

\medskip
\address
Aleksandar Ivi\'c, Katedra Matematike RGF-a
Universiteta u Beogradu, \DJ u\v sina 7, 11000 Beograd, Serbia
\endaddress
\email {\tt
aleksandar.ivic\@rgf.bg.ac.rs,  aivic\@matf.bg.ac.rs} \endemail
\dedicatory
\enddedicatory
\abstract
{Let $\D(x)$ denote the error term in the classical Dirichlet
divisor problem, and let the modified error term in the
divisor problem be $\D^*(x) =
 -\D(x)  + 2\D(2x) - \hf\D(4x)$. We show that
 $$
 \int_T^{T+H}\D^*\bigl(\frac{t}{2\pi}\bigr)|\zt|^2\d t \;\ll\; HT^{1/6}\log^{7/2}T
 \quad(T^{2/3+\e} \le H = H(T) \le T),
 $$
 $$
 \int_0^T\D(t)|\zt|^2\d t \;\ll\; T^{9/8}(\log T)^{5/2},
 $$
 and obtain  asymptotic formulae for
 $$
 \int_0^T{\Bigl(\D^*\bigl(\frac{t}{2\pi}\bigr)\Bigr)}^2|\zt|^2\d t,\quad
 \int_0^T{\Bigl(\D^*\bigl(\frac{t}{2\pi}\bigr)\Bigr)}^3|\zt|^2\d t.
 $$
 The importance of the $\Delta^*$-function comes from the fact that it is the analogue 
 of $E(T)$, the error term in the mean square formula for $|\zt|^2$. We also show, if
 $E^*(T) = E(T) - 2\pi \D^*(T/(2\pi))$, 
 $$
\int_0^T E^*(t)E^j(t)|\zt|^2\d t \;\ll_{j,\e}\; T^{7/6+j/4+\e}\quad(j= 1,2,3).
$$ }
\endabstract
\endtopmatter

\document
\head
1. Introduction
\endhead

As usual, let
$$
\D(x) \;:=\; \sum_{n\le x}d(n) - x(\log x + 2\gamma - 1)\qquad(x\geqslant 2)
\leqno(1.1)
$$
denote the error term in the classical Dirichlet divisor problem. Also let
$$
E(T) \;:=\;\int_0^T|\zt|^2\d t - T\Bigl(\log\bigl({T\over2\pi}\bigr) + 2\gamma - 1
\Bigr)\qquad(T\geqslant2)\leqno(1.2)
$$
denote the error term in the mean square formula for $|\zt|$.
Here $d(n)$ is the number of all positive divisors of
$n$, $\z(s)$ is the Riemann zeta-function, and $ \gamma = -\G'(1) = 0.577215\ldots\,$
is Euler's constant.  Long ago F.V. Atkinson [1] established a fundamental explicit formula
for $E(T)$ (see also [5, Chapter 15] and [7, Chapter 2]), which indicated a certain
 analogy between $\D(x)$ and $E(T)$. However, in this context it seems
that instead of the error-term function $\D(x)$ it
is more exact to work with the modified
function $\D^*(x)$ (see  M. Jutila [12], [13] and T. Meurman [14]), where
$$
\eqalign{
\D^*(x) :&= -\D(x)  + 2\D(2x) - \hf\D(4x)\cr&
= \hf\sum_{n\le4x}(-1)^nd(n) - x(\log x + 2\gamma - 1),\cr}
\leqno(1.3)
$$
since it turns out that $\D^*(x)$ is a better analogue of $E(T)$ than $\D(x)$.
Namely, M. Jutila (op. cit.) investigated both the
local and global behaviour of the difference
$$
E^*(t) \;:=\; E(t) - 2\pi\D^*\bigl({t\over2\pi}\bigr),\leqno(1.4)
$$
and in particular in [13] he proved that
$$
\int_T^{T+H}(E^*(t))^2\d t \;\ll_\e\; HT^{1/3}\log^3T+ T^{1+\e}\qquad(1\le H\le T).
\leqno(1.5)
$$
Here and later $\e$ denotes positive constants which are arbitrarily
small, but are not necessarily the same ones at each occurrence,
while $a(x) \ll_\e b(x)$ (same as $a(x) = O_\e(b(x)))$ means that
the $|a(x)| \le Cb(x)$ for some $C = C(\e) >0, x\ge x_0$.
The significance of (1.5) is that, in view of (see e.g., [5, Chapter 15])
$$
\int_0^T(\D^*(t))^2\d t \;\sim\; AT^{3/2},\quad
 \int_0^T E^2(t)\d t \;\sim \;BT^{3/2}\quad(A,B >0,\;
T\to\infty),\leqno(1.6)
$$
it transpires that $E^*(t)$ is in the mean square sense of a lower order of magnitude than
either $\D^*(t)$ or $E(t)$. We also refer the reader
to the review paper [18] of K.-M. Tsang on this subject.

\head
2. Statement of results
\endhead

Mean values (or moments) of $|\zt|$ represent  one of the central themes in the theory of $\z(s)$,
and they have been studied extensively.
There are two monographs dedicated solely to them: the author's [7], and that of K.
Ramachandra [17]. We are interested in obtaining mean value results for $\D^*(t)$ and $|\zt|^2$,
namely how the quantities in question relate to one another.
Our results are as follows.

\medskip
THEOREM 1. {\it For $T^{2/3+\e} \le H = H(T) \le T$ we have}
$$
 \int_T^{T+H}\D^*\bigl(\frac{t}{2\pi}\bigr)|\zt|^2\d t \;\ll\; HT^{1/6}\log^{7/2}T.\leqno(2.1)
 $$
\medskip
{\bf Remark 1}. If one uses the first formula in (1.6), the classical bound (see e.g., [5, Chapter 4])
$$
\int_0^T|\zt|^4\d t \;\ll\; T\log^4T \leqno(2.2)
$$ 
and the Cauchy-Schwarz inequality for integrals, one obtains
$$
\int_0^T\D^*\bigl(\frac{t}{2\pi}\bigr)|\zt|^2\d t \;\ll\; T^{5/4}\log^2T,
$$
which is considerably poorer than (2.1) for $H=T$, thus showing that this bound of Theorem 1 is non-trivial.

\medskip
THEOREM 2. {\it If $\gamma $ is Euler's constant and
$$
C := \frac{2\z^4(3/2)}{3\sqrt{2\pi}\z(3)} = 
\frac{2}{3\sqrt{2\pi}}\sum_{n=1}^\infty d^2(n)n^{-3/2}
= 10.3047\ldots\,, \leqno(2.3)
$$
then
$$
\int_0^T  {\Bigl(\D^*\bigl(\frac{t}{2\pi}\bigr)\Bigr)}^2|\zt|^2\d t =
\frac{C}{4\pi^2}T^{3/2}\Bigl(\log\frac{T}{2\pi} + 2\gamma - \frac{2}{3}\Bigr)
+ O_\e(T^{17/12+\e}).\leqno(2.4)
$$}

\medskip
{\bf Remark 2}. Note that (2.4) is a true asymptotic formula (17/12 = 3/2 - 1/12).
It would be interesting to analyze the error term in (2.4) and see how small it
can be, i.e., to obtain an omega-result (recall that $f(x) = \Omega(g(x))$ means
that $f(x) = o(g(x))$ does not hold as $x\to\infty$).

\medskip

THEOREM 3. {\it For some explicit constant $D > 0$ we have}
$$
\int_0^T  {\Bigl(\D^*\bigl(\frac{t}{2\pi}\bigr)\Bigr)}^3|\zt|^2\d t =
DT^{7/4}\Bigl(\log\frac{T}{2\pi} + 2\gamma - \frac{4}{7}\Bigr)
+ O_\e(T^{27/16+\e}).\leqno(2.5)
$$

\medskip
{\bf Remark 3}. Like (2.4), the formula in (2.5) is also
a true asymptotic formula (27/16 = 7/4 - 1/16).
Moreover, the main term is positive, which shows that, in the mean, 
$ \D^*\bigl(\frac{t}{2\pi}\bigr)$ is more biased towards positive values.

\medskip
In the most interesting case when $H =T$, Theorem 1 can be improved. Indeed, we have

\medskip
THEOREM 4. {\it We have
$$
\int_0^T\D(t)|\zt|^2\d t \;\ll\; T^{9/8}(\log T)^{5/2},\leqno(2.6)
$$
and} (2.6) {\it remains true if $\D(t)$ is replaced by $\D^*(t), \D(t/(2\pi))$ or 
$\D^*(t/(2\pi))$.}

\medskip
{\bf Remark 3}. The presence of $\D^*\bigl(\frac{t}{2\pi}\bigr)$ instead of the more
natural $\D^*(t)$ in (2.1), (2.4) and (2.5) comes from the defining relation (1.4). It would
be interesting to see what could be proved if in the integrals in (2.4) and (2.5) one
had $\D^*(t)$ (or $\D(t)$) instead of of $\D^*\bigl(\frac{t}{2\pi}\bigr)$. 

\medskip
{\bf Remark 4}. In the case of (2.1) (when $H=T$), Theorem 4 answers this question. However,
obtaining a short interval result for $\D(t)|\zt|^2$ is not easy. The method of proof
of Theorem 4 cannot be easily adapted to yield the analogues of (2.4) and (2.5) for
$\D(t)$ in place of $\D^*(t/(2\pi))$.

\medskip
There are some other integrals which may be bounded by the method used to prove previous theorems.
For example, one such result is

\medskip
THEOREM 5. {\it For $j = 1,2,3$ we have}
$$
\int_0^T E^*(t)E^j(t)|\zt|^2\d t \;\ll_{j,\e}\; T^{7/6+j/4+\e}.\leqno(2.7)
$$

\head
3. The necessary lemmas
\endhead
In this section we shall state some lemmas needed for the proof of our theorems.
The proofs of the theorems themselves will be given in Section 4.
\medskip
The first lemma embodies some bounds for the higher moments of $E^*(T)$. 

\medskip
LEMMA 1. {\it We have
$$
\int_0^T |E^*(t)|^3\d t \;\ll_\e\; T^{3/2+\e},\leqno(3.1)
$$
$$
\int_0^T |E^*(t)|^5\d t \;\ll_\e\; T^{2+\e},\leqno(3.2)
$$
and also}
$$
\int_0^T (E^*(t))^4\d t \;\ll_\e\; T^{7/4+\e}.\leqno(3.3)
$$

\medskip
The author proved (3.1) in [8, Part IV], and (3.2) in  [8, Part II]. The bound
(3.3) follows from (3.1) and (3.2) by the Cauchy-Schwarz inequality for integrals.

\medskip

For the mean square of $E(t)$ we need a more precise formula than (1.6). This is

\medskip
LEMMA 2. {\it With $C$ given by} (2.3) {\it we have}
$$
\int_0^TE^2(t)\d t = CT^{3/2} + R(T), \quad R(T) = O(T\log^4T).\leqno(3.4)
$$

\medskip
The first result on $R(T)$ is due to D.R. Heath-Brown [2], who obtained
$R(T) = O(T^{5/4}\log^2T)$. The sharpest known result at present is
$R(T) = O(T\log^4T)$, due independently to E. Preissmann [16]
and the author [7, Chapter 2]. 

\medskip
For the mean square of $E^*(t)$ we have a result which is different from (3.4). This is

\medskip
LEMMA 3. {\it We have
$$
\int_0^T (E^*(t))^2\d t \;=\; T^{4/3}P_3(\log T) + O_\e(T^{7/6+\e}),\leqno(3.5)
$$
where $P_3(y)$ is a polynomial of degree three in $y$ with positive leading 
coefficient, and all its coefficients may be evaluated explicitly.}

\medskip
This formula was proved by the author in [9]. It sharpens (1.4) when $H=T$.
It seems likely that the error term in (3.5) is $O_\e(T^{1+\e})$, but this 
seems difficult to prove.

\medskip
LEMMA 4. {\it We have
$$
\int_0^T|\zt|^4 \d t = TQ_4(\log T) + O(T^{2/3}\log^8T),\leqno(3.6)
$$
where $Q_4(x)$ is an explicit polynomial of degree four in $x$ with
leading coefficient $1/(2\pi^2)$.
}

\medskip
This result was proved first (with error term $O(T^{2/3}\log^CT)$) by
Y. Motohashi and the author [10]. The value $C=8$ was given by Y. Motohashi 
in his monograph [15].

\medskip
LEMMA 5. {\it For $1 \le N \ll x$ we have}
$$
\D^*(x) = {1\over\pi\sqrt{2}}x^{1\over4}
\sum_{n\le N}(-1)^nd(n)n^{-{3\over4}}
\cos(4\pi\sqrt{nx} - {\txt{1\over4}}\pi) +
O_\e(x^{{1\over2}+\e}N^{-{1\over2}}).
\leqno(3.7)
$$
\medskip
The expression for $\D^*(x)$ (see [5, Chapter 15]) is the analogue
of the classical truncated Vorono{\"\i} formula for $\D(x)$ (ibid. Chapter 3),
which is the expression in (3.7) without $(-1)^n$.

\medskip

LEMMA 6. {\it We have
$$
\int_0^TE(t)|\zt|^2\d t = \pi T\Bigl(\log\frac{T}{2\pi} + 2\gamma-1\Bigr) + U(T),
\leqno(3.8)
$$
where
$$
U(T) = O(T^{3/4}\log T),\quad U(T) = \Omega_\pm(T^{3/4}\log T).
$$}

\medskip
The asymptotic formula (3.8) is due to the author [6]. Here the symbol $f(x) = \Omega_\pm(g(x))$
has its standard meaning, namely that both $\limsup_{x\to\infty}f(x)/g(x) > 0$ and 
$\liminf_{x\to\infty}f(x)/g(x)<0$ holds.

\medskip
LEMMA 7. {\it We have
$$
\eqalign{&
\int_1^TE^3(t)\d t = C_1T^{7/4} + O_\e(T^{5/3+\e}),
\cr &
\int_1^TE^4(t)\d t = C_2T^2 + O_\e(T^{23/12+\e}),\cr}\leqno(3.9)
$$
where $C_1, C_2$ are certain explicit, positive constants.}

\medskip
These asymptotic formulae are due to P. Sargos and the author [11].

\medskip
LEMMA 8. {\it We have
$$
\sum_{n\le x}d^2(n) \;=\; \frac{1}{\pi^2}x\log^3x + O(x\log^2x).\leqno(3.10)
$$
}
\medskip
This is a well-known elementary formula; see e.g., page 141 of [5].

\medskip
LEMMA 9. {\it For real $k \in [0,9]$ the limits
$$
E_k 
\;:=\; \lim_{T\to\infty}T^{-1-k/4}\int_0^T|E(t)|^k\d t
$$
exist.}

\medskip
This is a result of D.R. Heath-Brown [4]. The limits of moments without
absolute values also exist when  $k=1,3,5,7$ or 9.

\medskip
LEMMA 10. {\it For $4\le A \le 12$ we have
$$
\int_0^T |\zt|^A\d t \;\ll_A\; T^{1 + \frac{1}{8}(A-4)}\log^{C(A)}T\leqno(3.11)
$$
with some positive constant $C(A)$.}

\medskip
These are at present the strongest upper bounds  for moments of $|\zt|$ for
the range in question. They follow by convexity from the fourth moment bound (2.2)
and the twelfth moment
$$
\int_0^T|\zt|^{12}\d t \;\ll\; T^2\log^{17}T
$$
of D.R. Heath-Brown [3] (see e.g., [5, Chapter 8] for more details).

\head
4. Proofs of the Theorems
\endhead

We begin with the proof of (2.1). We start from
$$
\eqalign{
\int_T^{T+H}E^*(t)|\zt|^2\d t &\ll \left\{\int_T^{T+H}(E^*(t))^2\d t
\int_T^{T+H}|\zt|^4\d t\right\}^{1/2}\cr&
\ll \left(HT^{1/3}\log^3T\cdot H\log^4T\right)^{1/2} = HT^{1/6}\log^{7/2}T.
\cr}\leqno(4.1)
$$
Here we assumed that $T^{2/3+\e} \le H = H(T) \le T$ and used   
(3.6) of Lemma 4, (1.5) and the Cauchy-Schwarz inequality for
integrals. On the other hand, by the defining relation (1.4) we have
$$
\eqalign{
\int_T^{T+H} E^*(t)|\zt|^2\d t &= \int_T^{T+H} E(t)|\zt|^2\d t \cr&- 
2\pi\int_T^{T+H}\D^*\bigl(\frac{t}{2\pi}\bigr)|\zt|^2\d t.\cr} \leqno(4.2)
$$ 
Using (3.8) of Lemma 6 and (4.1), we obtain then from (4.2)
$$
\eqalign{&
2\pi\int_T^{T+H}\D^*\bigl(\frac{t}{2\pi}\bigr)|\zt|^2\d t = \int_T^{T+H} E^*(t)|\zt|^2\d t 
\cr&- \int_T^{T+H} E(t)|\zt|^2\d t\cr&
= O(HT^{1/6}\log^{7/2}T) + \pi t\Bigl(\log\frac{t}{2\pi} + 2\gamma-1\Bigr)\Bigl|_T^{T+H} + O(T^{3/4}\log T)
\cr&
= O(HT^{1/6}\log^{7/2}T) + O(H\log T)+ O(T^{3/4}\log T)\cr&
\ll HT^{1/6}\log^{7/2}T,
\cr}
$$
since $T^{2/3+\e} \le H = H(T) \le T$. This completes the proof of Theorem 1.

\medskip
The proof of Theorem 2 is somewhat more involved. It suffices to consider the integral from $T$ to
$2T$, and then at the end of the proof to replace $T$ by $T2^{-j}$ and sum the resulting expressions
when $j = 1,2,\ldots\;$. First, by squaring (1.4), we have
$$
\eqalign{&
\int_T^{2T}(E^*(t))^2|\zt|^2\d t = \int_T^{2T}(E(t))^2|\zt|^2\d t\cr&
- 2\int_T^{2T}E(t)2\pi\D^*\bigl(\frac{t}{2\pi}\bigr)|\zt|^2\d t
+4\pi^2\int_T^{2T}\Bigl(\D^*\bigl(\frac{t}{2\pi}\bigr)\Bigr)^2|\zt|^2\d t.
\cr}\leqno(4.3)
$$
The expression in the  middle of the right-hand side of (4.3) equals, on differentiating (1.2),
$$
\eqalign{&
-2\int_T^{2T}E(t)2\pi\D^*\bigl(\frac{t}{2\pi}\bigr)\Bigl(\log\frac{t}{2\pi}+2\gamma + E'(t)\Bigr)\d t
\cr&
= -2\int_T^{2T}E(t)2\pi\D^*\bigl(\frac{t}{2\pi}\bigr)\Bigl(\log\frac{t}{2\pi}+2\gamma\Bigr)\d t + J(T),
\cr}\leqno(4.4)
$$
say, where
$$
J(T):= -2\int_T^{2T}E(t)2\pi\D^*\bigl(\frac{t}{2\pi}\bigr)E'(t)\d t.\leqno(4.5)
$$
To bound $J(T)$ we use Lemma 5 with $N = N(T), 1\ll N \ll T$, where $N$ will be determined
a little later. The error term in (3.7) trivially makes a contribution which is
$$
\ll_\e \int_T^{2T}|E(t)|\Bigl(|\zt|^2 + \log T\Bigr)T^{1/2+\e}N^{-1/2} 
\ll_\e T^{7/4+\e}N^{-1/2}\leqno(4.6)
$$
on using the second formula in (1.6), (2.2) and the Cauchy-Schwarz inequality for integrals.
There remains the contribution of a multiple of
$$
\int_T^{2T}{\bigl(E^2(t)\bigr)}^{'}t^{1/4}\sum_{n\le N}(-1)^nd(n)n^{-3/4}\cos(\sqrt{8\pi nt}-\pi/4)\d t.
$$
This is integrated by parts. The integrated terms are $\ll T^{11/12}N^{1/3}\log T$, by using the standard
estimate $E(T) \ll T^{1/3}$ (see e.g., [5, Chapter 15]) and trivial estimation. 
The main contribution comes from the
differentiation of the sum over $n$. Its contribution will be, with $n\sim K$ meaning
that $K < n \le K' \le2K$,
$$
\eqalign{&
\ll T^{-1/4}\int_T^{2T}E^2(t)\Bigl|\sum_{n\le N}(-1)^n d(n)n^{-1/4}\exp(i\sqrt{8\pi nt})\Bigr|\d t
\cr&
\le T^{-1/4}\left\{\int_T^{2T}E^4(t)\d t\int_T^{2T}\Bigl|
\sum_{n\le N}(-1)^n d(n)n^{-1/4}\exp(i\sqrt{8\pi nt})\Bigr|^2\d t\right\}^{1/2}\cr&
\ll T^{3/4}\left(\int_T^{2T}\sum_{n\le N}d^2(n)n^{-1/2}\d t + 
\log^2 T\max_{K\ll N}\sum_{m\not= n\sim K}K^{\e-1/2} \frac{\sqrt{T}}{|\sqrt{m}-\sqrt{n}|}\right)^{1/2}
\cr&
\ll T^{3/4}(TN^{1/2}\log^3T + T^{1/2+\e}N)^{1/2} \ll T^{5/4}N^{1/4}\log^{3/2}T
\cr}\leqno(4.7)
$$
for $T^\e \le N = N(T) \le T^{1-\e}$. Here we used the standard first derivative test (see e.g., Lemma 2.1
of [6]) for exponential integrals, Lemma 7, (3.10) and
$$
\sum_{m\not= n\sim K} \frac{1}{|\sqrt{m}-\sqrt{n}|}
\;\ll\; \sum_{n\sim K}\sum_{m\sim K, m\ne n}\frac{\sqrt{K}}{|m -n|} \;\ll\; K^{3/2}\log K.
$$ 
From (4.6) and (4.7) we see that the right choice for $N$ should be if we have
$$
T^{7/4}N^{-1/2} \;=\; T^{5/4}N^{1/4},\quad N \;=\; T^{2/3},
$$
and with this choice of $N$ we obtain $T^{11/12}N^{1/3} = T^{41/36}\;(41/36 < 17/12)$, and
$$
J(T) \;\ll_\e\; T^{17/12+\e}. \leqno(4.8)
$$
In view of (1.4), the formula (4.3) and the bound (4.8) give
$$
\eqalign{&
4\pi^2\int_T^{2T}{\Bigl(\D^*\bigl(\frac{t}{2\pi}\bigr)\Bigr)}^2|\zt|^2\d t = O_\e(T^{17/12+\e})\cr&
+ 2\int_T^{2T}E(t)2\pi\D^*\bigl(\frac{t}{2\pi}\bigr)\Bigl(\log\frac{t}{2\pi}+2\gamma\Bigr)\d t
- \int_T^{2T}E^2(t)|\zt|^2\d t \cr&
+ \int_T^{2T}{(E^*(t))}^2|\zt|^2\d t = O_\e(T^{17/12+\e}) +2I_1 - I_2 + I_3, 
\cr}\leqno(4.9)
$$
say. On using (3.3) of Lemma 1, (2.2) and the Cauchy-Schwarz inequality we obtain
$$
I_3 \ll \left\{\int_T^{2T}(E^*(t))^4\d t\int_T^{2T}|\zt|^4\d t \right\}^{1/2} \ll_\e T^{11/8+\e}.
$$
Further we have
$$
\eqalign{&
2I_1 - I_2 = 2\int_T^{2T}E(t)\Bigl(E(t)- E^*(t)\Bigr)
\Bigl(\log\frac{t}{2\pi}+2\gamma\Bigr)\d t - I_2
\cr&
= \int_T^{2T}E^2(t)\left\{2\Bigl(\log\frac{t}{2\pi}+2\gamma\Bigr)-|\zt|^2\right\}\d t
\cr&
-2\int_T^{2T}E(t)E^*(t)\Bigl(\log\frac{t}{2\pi}+2\gamma\Bigr)\d t\cr&
= \int_T^{2T}E^2(t)\Bigl(\log\frac{t}{2\pi}+2\gamma- E'(t)\Bigr)\d t
- 2\int_T^{2T}E(t)E^*(t)\Bigl(\log\frac{t}{2\pi}+2\gamma\Bigr)\d t.
\cr}
$$
The last integral is, by Lemma 2, Lemma 3 and the Cauchy-Schwarz inequality for integrals,
$$
\ll \log T{\left\{\int_T^{2T}E^2(t)\d t\int_T^{2T}(E^*(t))^2\d t\right\}}^{1/2} 
\ll T^{17/12}\log^{5/2}T.
$$
On the other hand,
$$
\eqalign{&
\int_T^{2T}E^2(t)\Bigl(\log\frac{t}{2\pi}+2\gamma- E'(t)\Bigr)\d t
\cr&
= \int_T^{2T}E^2(t)\Bigl(\log\frac{t}{2\pi}+2\gamma\Bigr)\d t - \txt{\frac{1}{3}}E^3(t)\Bigl|_T^{2T}
\cr&
= \int_T^{2T}E^2(t)\Bigl(\log\frac{t}{2\pi}+2\gamma\Bigr)\d t + O(T).
\cr}\leqno(4.10)
$$
To evaluate the last integral in (4.10) we use Lemma 6 and integration by parts. 
This shows that the integral in question is
$$
\eqalign{&
\bigl(Ct^{3/2} + R(t)\bigr)\Bigl(\log\frac{t}{2\pi}+2\gamma\Bigr)\Bigl|_T^{2T}
- \int_T^{2T}\Bigl(Ct^{1/2} + \frac{R(t)}{t}\Bigr)\d t
\cr&
= Ct^{3/2}\Bigl(\log\frac{t}{2\pi}+2\gamma\Bigr)\Bigl|_T^{2T} + O(T\log^5T) 
-\txt{\frac{2}{3}}Ct^{3/2}\Bigl|_T^{2T}
\cr&
= Ct^{3/2}\Bigl(\log\frac{t}{2\pi}+2\gamma- {\frac{2}{3}}\Bigr)\Bigl|_T^{2T}
\,+\, O(T\log^5T). 
\cr}
$$
It transpires from (4.9) and (4.10) that
$$
4\pi^2\int_T^{2T}{\bigl(\D^*\bigl(\frac{t}{2\pi}\bigr)\bigr)}^2|\zt|^2\d t =
Ct^{3/2}\Bigl(\log\frac{t}{2\pi}+2\gamma- {\frac{2}{3}}\Bigr)\Bigl|_T^{2T}
+ O_\e(T^{17/12+\e}),
$$
which gives at once (2.4) of Theorem 2. 

\medskip
We turn now to the proof of Theorem 3. The basic idea is analogous to the one used in the
proof of Theorem 2, so that we shall be relatively brief. The integral in (2.5)
equals $1/(8\pi^3)$ times
$$
\int_0^T\Bigl\{E^3(t)- 3E^*(t)E^2(t) + 3(E^*(t))^2E(t)- (E^*(t))^3 \Bigr\}|\zt|^2\d t.\leqno(4.11)
$$
The main term in (2.5) comes from
$$
\eqalign{&
\int_0^T E^3(t)|\zt|^2\d t = \int_0^T E^3(t)\Bigl(\log\frac{t}{2\pi}+2\gamma-E'(t)\Bigr)\d t
\cr&
= C_1T^{7/4}\Bigl(\log\frac{T}{2\pi}+2\gamma\Bigr) - \int_1^T C_1t^{3/4}\d t + O_\e(T^{5/3+\e})
\cr&
= C_1T^{7/4}\Bigl(\log\frac{T}{2\pi}+2\gamma-\frac{4}{7}\Bigr) + O_\e(T^{5/3+\e}),
\cr}
$$
where (3.9) of Lemma 7 was used. By H\"older's inequality for integrals, (3.3) of Lemma 1 and
(3.10) of Lemma 3 (with $A=5$) we obtain
$$
\eqalign{&
\int_0^T(E^*(t))^3|\zt|^2\d t \ll \left(\int_0^T|E^*(t)|^5\d t\right)^{3/5}
\left(\int_0^T|\zt|^5\d t\right)^{2/5}
\cr&
\ll_\e\; T^{6/5+ 9/20+\e} = T^{33/20+\e}.
\cr}
$$
Similarly we obtain
$$
\eqalign{&
\int_0^T(E^*(t))^2E(t)|\zt|^2\d t 
\cr&\ll \left(\int_0^T|E^*(t)|^{16/3}\d t\right)^{3/8}\left(\int_0^T E^8(t)\d t\right)^{1/8}
\left(\int_0^T|\zt|^4\d t\right)^{1/2}
\cr&
\ll_\e\; T^{\frac{3}{8}(2+\frac{1}{9})+\frac{3}{8}+\frac{1}{2}+\e} = T^{\frac{5}{3}+\e},
\cr}
$$
where we used (3.2) and the fact that $E^*(T) \ll T^{1/3}$, which follows from the definition of
$E^*$ and the classical estimates $\D(x) \ll x^{1/3}, E(T) \ll T^{1/3}$.
Finally, by using (3.3), Lemma 9 with $k=8$ and (2.2), we obtain
$$
\eqalign{&
\int_0^TE^*(t)E^2(t)|\zt|^2\d t 
\cr&
\ll \left(\int_0^T|E^*(t)|^4\d t\right)^{1/4}
\left(\int_0^T E^8(t)\d t\right)^{1/4}
\left(\int_0^T|\zt|^4\d t\right)^{1/2}
\cr&
\ll_\e T^{7/16+3/4+ 1/2+\e} = T^{27/16+\e}.
\cr}
$$
Since $27/16 = 1.6875 > 5/3 > 33/20 = 1.65$, we obtain easily the assertion of Theorem 3.

\medskip
We shall  prove now (2.6) of Theorem 4. We suppose $T \le t \le 2T$ and
take $N = T$ in (3.7) of Lemma 5. This holds both for $\D^*(x)$ and $\D(x)$, and
one can see easily that the proof remains valid if we have an additional factor
of $1/(2\pi)$ in the argument of $\D^*$ or $\D$ (or any constant $c>0$, for that matter). 
Thus we start from
$$
\eqalign{
\D(t) &= \frac{t^{1/4}}{\pi\sqrt{2}}\sum_{n\le T}d(n)n^{-3/4}\cos(4\pi\sqrt{nt}-\pi/4) + O_\e(T^\e)
\cr&
= \frac{t^{1/4}}{\pi\sqrt{2}}\left(\sum_{n\le G}\cdots + \sum_{G < n\le T}\cdots\right) + O_\e(T^\e),
\cr}\leqno(4.12)
$$
say, where $T^\e \le G = G(T) \le T^{1-\e}$, and $G$ will be determined a little later.
The error term in (4.12) makes a contribution of $O_\e(T^{1+\e})$ to (2.6). We have
$$
\eqalign{&
\int_T^{2T} t^{1/4}\sum_{n\le G}d(n)n^{-3/4}\cos(4\pi\sqrt{nt}-\pi/4)|\zt|^2\d t\cr& 
= \int_T^{2T}t^{1/4}\Bigl(\log{\frac{t}{2\pi}}
 + 2\gamma + E'(t)\Bigr)\sum_{n\le G}d(n)n^{-3/4}\cos(4\pi\sqrt{nt}-\pi/4)
\d t\cr&
= I_1 + I_2,
\cr}\leqno(4.13)
$$
say. By the first derivative test
$$
\eqalign{
I_1 :=&
 \int_T^{2T}t^{1/4}\Bigl(\log{\frac{t}{2\pi}}
 + 2\gamma \Bigr)\sum_{n\le G}d(n)n^{-3/4}\cos(4\pi\sqrt{nt}-\pi/4)\d t\cr&
\ll T^{1/4}\log T\cdot \sum_{n\le G}d(n)n^{-3/4}T^{1/2}n^{-1/2} \ll T^{3/4}\log T,
\cr}
$$
since $\sum_{n\ge1}d(n)n^{-\a}$ converges for $\a>1$. The integral $I_2$, namely
$$
I_2:= \int_T^{2T}t^{1/4}E'(t)\sum_{n\le G}d(n)n^{-3/4}\cos(4\pi\sqrt{nt}-\pi/4)
\d t
$$
is integrated by parts. The integrated terms are trivially $O(T)$, and there remains
$$
\eqalign{
&
 -\int_T^{2T}\frac{1}{4}t^{-3/4}E(t)\sum_{n\le G}d(n)n^{-3/4}\cos(4\pi\sqrt{nt}-\pi/4)\d t\cr&
+2\pi \int_T^{2T}t^{-1/4}E(t)\sum_{n\le G}d(n)n^{-1/4}\sin(4\pi\sqrt{nt}-\pi/4)\d t.
\cr}\leqno(4.13)
$$
Both integrals in (4.13) are estimated analogously, and clearly it is the latter which is larger.
By the Cauchy-Schwarz inequality for integrals it is
$$
\ll\; T^{-1/4}{(J_1J_2)}^{1/2},
$$
where
$$
\eqalign{
J_1 &:=
 \int_T^{2T}{\Bigl|\sum_{n\le G}d(n)n^{-1/4}\E^{4\pi i\sqrt{nt}}\Bigr|}^2\d t\cr
J_2 &:= \int_T^{2T}E^2(t)\d t \;\ll\; T^{3/2},
\cr}
$$
on using Lemma 2 in bounding $J_2$. Using the first derivative test and (3.10) of Lemma 8, we find that
$$
\eqalign{
J_1 &=
T\sum_{n\le G}d^2(n)n^{-1/2} + \sum_{m\ne n\le G}\frac{d(m)d(n)}{(mn)^{1/4}}
 \int_T^{2T}\E^{4\pi i(\sqrt{m}-\sqrt{n})\sqrt{t}}\d t\cr&
\ll TG^{1/2}\log^3T + T^{1/2}\sum_{m\ne n\le G}\frac{d(m)d(n)}{(mn)^{1/4}|\sqrt{m}-\sqrt{n}|}.
\cr}
$$
When $n/2 < m \le 2n$ the contribution of the last double sum is
$$
\ll_\e T^{1/2}\sum_{n\le G}n^{\e-1/2}n^{1/2}\sum_{n/2<m\le 2n, m\ne n}\frac{1}{|m-n|} \ll_\e T^{1/2+\e}G.
$$
If $m \le n/2$ then $|\sqrt{m}-\sqrt{n}|^{-1} \ll n^{-1/2}$, and when $m>2n$ it is $\ll m^{-1/2}$. Thus the
total contribution of the double sum above is certainly
$$
\ll_\e T^{1/2+\e}G \ll TG^{1/2}\log^3T\qquad(T^\e \le G = G(T) \le T^{1-\e}).
$$
We infer that
$$
T^{-1/4}{(J_1J_2)}^{1/2} \ll T^{-1/4}(TG^{1/2}\log^3T\cdot T^{3/2})^{1/2} = TG^{1/4}(\log T)^{3/2}.
$$
In a similar vein it is found that
$$
\eqalign{&
\int_T^{2T} t^{1/4}\sum_{G<n\le T}d(n)n^{-3/4}\cos(4\pi\sqrt{nt}-\pi/4)|\zt|^2\d t\cr& 
\ll T^{1/4}\left\{\int_T^{2T}\Bigl|\sum_{G<n\le T}d(n)n^{-3/4}\E^{4\pi i\sqrt{nt}}\Bigr|^2\d t
\int_T^{2T}|\zt|^4\d t\right\}^{1/2}\cr&
\ll T^{3/4}\log^2T\left\{\int\limits_T^{2T}\left(\sum_{n>G}\frac{d^2(n)}{n^{3/2}} +
\sum_{G < m\ne n\le T}\frac{d(m)d(n)}{(mn)^{3/4}}\E^{4\pi i(\sqrt{m}-\sqrt{n})\sqrt{t})}\right)\d t\right\}^{1/2}\cr&
\ll_\e T^{3/4}\log^2T\Bigl\{TG^{-1/2}\log^3T + T^{1/2+\e}\Bigr\}^{1/2} \ll T^{5/4}G^{-1/4}(\log T)^{7/2}.
\cr}
$$
We finally infer that, for $T^\e \le G = G(T) \le T^{1-\e}$,
$$
\int\limits_T^{2T}\D(t)|\zt|^2\d t \ll TG^{1/4}(\log T)^{3/2} + T^{5/4}G^{-1/4}(\log T)^{7/2} \ll T^{9/8}(\log T)^{5/2}
$$
with the choice $G = T^{1/2}\log^4T$. This leads to (2.6) on replacing $T$ by $T2^{-j}$ and adding the resulting
estimates.

\medskip
{\bf Corollary}. We have
$$
\int_0^TE^*(t)|\zt|^2\d t \;\ll\; T^{9/8}(\log T)^{5/2}.\leqno(4.14)
$$
Namely
$$
\int_T^{2T}E^*(t)|\zt|^2\d t = \int_T^{2T}\Bigl\{E(t) - 2\pi\D^*(t/(2\pi))\Bigr\}|\zt|^2\d t.
$$
The integral with $\D^*$ is $\ll T^{9/8}(\log T)^{5/2}$ by Theorem 4. There remains
$$
\int_T^{2T}E(t)|\zt|^2\d t = \pi T\Bigl(\log\frac{2T}{\pi} + 2\gamma-1\Bigr) + O(T^{3/4}\log T)
= O(T\log T)
 $$
 by (3.8) of Lemma 6. This gives
 $$
\int_T^{2T}E^*(t)|\zt|^2\d t \;\ll\; T^{9/8}(\log T)^{5/2}.
$$
To complete the proof of (4.14), again 
 one replaces $T$ by $T2^{-j}$ and adds the resulting
estimates.

\medskip
It remains to prove (2.7) of Theorem 5 (the bound (4.14) gives a
result when $j=0$). 
The proof  is analogous to the proofs given before, so we shall be brief. We have
$$
\eqalign{&
\int_T^{2T} E^*(t)E^j(t)|\zt|^2\d t\cr& 
= \int_T^{2T} E^*(t)E^j(t)\Bigl(\log{\frac{t}{2\pi}} + 2\gamma + E'(t)\Bigr)\d t = I' + I'',
\cr}
$$
say. By the Cauchy-Schwarz inequality for integrals, Lemma 3 and Lemma 9 (with $k=2j$), it follows that
$$
\eqalign{I'&
:= \int_T^{2T} E^*(t)E^j(t)\Bigl(\log{\frac{t}{2\pi}} + 2\gamma\Bigr)\d t\cr& 
\ll \log T{\Bigl\{\int_T^{2T}(E^*(t))^2\d t \int_T^{2T}E^{2j}(t)\d t\Bigr\}}^{1/2}\cr&
\ll \log T{(T^{4/3}\log^3T\cdot T^{1+j/2})}^{1/2} = T^{7/6+j/4}\log^{5/2}T.
\cr}
$$
On the other hand, by (1.4) we have
$$
\eqalign{I''&
:= \int_T^{2T} E^*(t)E^j(t)E'(t)\d t\cr& 
= \int_T^{2T} E^{j+1}(t)E'(t)\d t - 2\pi \int_T^{2T} \D^*\bigl( \frac{t}{2\pi}\bigr)E^j(t)E'(t)\d t.
\cr}\leqno(4.15)
$$
Note that 
$$
\int_T^{2T} E^{j+1}(t)E'(t)\d t = \textstyle\frac{1}{j+2}E^{j+2}(t)\Bigl|_T^{2T} = O(T^{(j+2)/3)}),
$$
and $(j+2)/3 < 7/6+j/4$ for $0 < j \le 6$.
By using (3.7) of Lemma 4 it is seen that the last integral in (4.15) is a multiple of
$$
\int_T^{2T}t^{1/4}\sum_{n\le N}(-1)^nd(n)n^{-3/4}\cos(\sqrt{8\pi nt}-\pi/4)E^j(t)E'(t)\d t + 
{\Cal J}(T),\leqno(4.16)
$$
say, where $T^\e \le N = N(T) \le T^{1-\e}$. Using Lemma 9 we have
$$
\eqalign{
{\Cal J}(T) &\ll_\e T^{1/2+\e}N^{-1/2}\int_T^{2T} |E^j(t)||E'(t)|\d t \cr&
\ll_\e T^{1/2+\e}N^{-1/2}{\left\{\int_T^{2T}E^{2j}(t)\d t
\int_T^{2T}\Bigl(\log^2T+|\zt|^4\Bigr)\d t\right\}}^{1/2}
\cr&
\ll_\e T^{1/2+\e}N^{-1/2}{\bigl(T^{1+j/2}\cdot T\log^4T\bigr)}^{1/2} = T^{3/2+ j/4+\e}N^{-1/2}.
\cr}
$$

The remaining integral in (4.16) is again integrated by parts. 
The major contribution will come from a multiple of
$$
\eqalign{&
 \int_T^{2T} E^{j+1}(t)t^{-1/4}\sum_{n\le N}(-1)^nd(n)n^{-1/4}\sin(\sqrt{8\pi nt}-\pi/4) \d t\cr& 
\ll T^{-1/4}{\Bigl\{\int_T^{2T}E^{2j+2}(t)\d t
\int_T^{2T}\Bigl|\sum_{n\le N}(-1)^nd(n)n^{-1/4}\E^{i\sqrt{8\pi nt}}\Bigr|^2\d t\Bigr\}}^{1/2}\cr&
\ll T^{-1/4}\bigl\{T^{1+(j+1)/2}\cdot TN^{1/2}\log^3T\bigr\}^{1/2} = T^{3/4+(j+1)/4}N^{1/4}\log^{3/2}T,
\cr}
$$
where Lemma 9 was used with $k=2j+2\le 8$. The choice $N = T^{2/3}$ gives
$$
T^{3/4+(j+1)/4}N^{1/4} = T^{3/2+ j/4+\e}N^{-1/2} = T^{7/6+j/4},
$$
as asserted by Theorem 5.
The bound in (2.7) is an expected one, since (in the mean square sense) $E^*(t)$ is
of the order $\ll t^{1/6}\log^{3/2}t$, $E(t)$ is of the order $\ll t^{1/4}$, and $|\zt|^2$ is
of logarithmic order. However, by H\"older's inequality for integrals (2.7) does not follows directly,
since it would require the (yet unknown) higher moments of $|\zt|$.

\medskip

\vfill
 \eject
\medskip
\Refs
\medskip

\item{[1]} F.V. Atkinson, {\it The mean value of the Riemann zeta-function},
Acta Math. {\bf81}(1949), 353-376.

\item{[2]} D.R. Heath-Brown, {\it The mean value theorem 
for the Riemann zeta-function}, Mathematika {\bf25}(1978), 177-184.

\item{[3]} D.R. Heath-Brown, {\it The twelfth power moment of the Riemann
zeta-function}, Quart. J. Math. (Oxford) {\bf29}(1978), 443-462.
   
\item{[4]} D.R. Heath-Brown, 
{\it The distribution and moments of the error term in the Dirichlet 
divisor problems}, Acta Arith. {\bf60}(1992), 389-415.

\item{[5]} A. Ivi\'c, {\it The Riemann zeta-function}, John Wiley \&
Sons, New York, 1985 (2nd ed. Dover, Mineola, New York, 2003).

\item{[6]} A. Ivi\'c, {\it On some integrals involving the mean square 
formula for the
    Riemann zeta-function}, Publications Inst. Math. (Belgrade)
    {\bf46(60)} (1989), 33-42.
    
\item {[7]} A. Ivi\'c,  {\it Mean values of the Riemann zeta-function},
LN's {\bf 82},  Tata Inst. of Fundamental Research,
Bombay,  1991 (distr. by Springer Verlag, Berlin etc.).

\item{[8]} A. Ivi\'c, {\it On the Riemann zeta-function and the divisor problem},
Central European J. Math. {\bf(2)(4)}\ (2004), 1-15,  II, ibid.
{\bf(3)(2)}(2005), 203-214,  III, Annales Univ.
Sci. Budapest, Sect. Comp. {\bf29}(2008), 3-23,
and IV, Uniform Distribution Theory {\bf1}(2006), 125-135.

\item{[9]} A. Ivi\'c, {\it On the mean square of the zeta-function and
the divisor problem}, Annales  Acad. Scien. Fennicae Mathematica {\bf23}(2007), 1-9.

\item{[10]} A. Ivi\'c and Y. Motohashi, {\it On the fourth power moment of the
Riemann zeta-function}, J. Number Theory {\bf51}(1995), 16-45.

\item{[11]} A. Ivi\'c and P. Sargos, {\it On the higher power moments of the
 error term in the divisor problem}, Illinois J. Math. {\bf81}(2007), 353-377.

\item{[12]} M. Jutila, {\it Riemann's zeta-function and the divisor problem},
Arkiv Mat. {\bf21}(1983), 75-96 and II, ibid. {\bf31}(1993), 61-70.

\item{[13]} M. Jutila, {\it On a formula of Atkinson},
Topics in classical number theory, Colloq. Budapest 1981,
Vol. I, Colloq. Math. Soc. J\'anos Bolyai {\bf34}(1984), 807-823.


\item{[14]} T. Meurman, {\it A generalization of Atkinson's formula to
$L$-functions}, Acta Arith. {\bf47}(1986), 351-370.

\item{[15]} Y. Motohashi, {\it Spectral theory of the Riemann zeta-function}, Cambridge
University Press, Cambridge, 1997.

\item{[16]} E. Preissmann, {\it Sur la moyenne quadratique du terme de 
reste du probl\'eme du cercle},
C. R. Acad. Sci. Paris S\'er. I Math. {\bf306}(1988), no. 4, 151-154.

\item{[17]} K.  Ramachandra, {\it On the mean-value and omega-theorems
for the Riemann zeta-function}, LN's {\bf85}, Tata Inst. of Fundamental Research
(distr. by Springer Verlag, Berlin etc.), Bombay, 1995.

\item{[18]} K.-M. Tsang, {\it Recent progress on the Dirichlet divisor problem and the 
mean square of the Riemann zeta-function}, Sci. China Math. {\bf53}(2010), no. 9, 2561-2572.

\vfill

\endRefs

\enddocument

\bye